\documentclass[12pt]{article}
\usepackage{amssymb, amsmath, latexsym, a4}
\usepackage[latin1]{inputenc}
\usepackage{graphicx}
\usepackage[all]{xy}

\newcommand{\rdg}{\hfill $\Box $}

\newtheorem{De}{Definition}[section]
\newtheorem{Th}[De]{Theorem}
\newtheorem{Pro}[De]{Proposition}
\newtheorem{Le}[De]{Lemma}
\newtheorem{Co}[De]{Corollary}
\newtheorem{Rem}[De]{Remark}
\newtheorem{Ex}[De]{Example}

\newcommand{\comp}{\raisebox{0.2mm}{\ensuremath{\scriptstyle{\circ}}}}

\renewcommand{\Im}{\ensuremath{\mathrm{Im}}}

\newcommand{\A}{\ensuremath{\mathcal{A}}}

\newcommand{\K}{\ensuremath{\mathbb{K}}}
\newcommand{\B}{\ensuremath{\mathcal{B}}}

\newcommand{\Lie}{\ensuremath{\mathsf{Lie}}}

\newcommand{\Lieh}{\ensuremath{\mathfrak{h}}}
\newcommand{\Lieg}{\ensuremath{\mathfrak{g}}}
\newcommand{\Lieq}{\ensuremath{\mathfrak{q}}}
\newcommand{\Liea}{\ensuremath{\mathfrak{a}}}
\newcommand{\Liem}{\ensuremath{\mathfrak{m}}}
\newcommand{\Lien}{\ensuremath{\mathfrak{n}}}
\newcommand{\Lief}{\ensuremath{\mathfrak{f}}}
\newcommand{\Lies}{\ensuremath{\mathfrak{s}}}
\newcommand{\Lier}{\ensuremath{\mathfrak{r}}}

\newcommand{\Leib}{\ensuremath{\mathsf{Leib}}}

\newcommand{\ze}{\cal Z}

\newbox\pullbackbox
\setbox\pullbackbox=\hbox{\xy 0;<1mm,0mm>: \POS(4,0)\ar@{-} (0,0) \ar@{-} (4,4)
\endxy}

\newdir{>>}{{}*!/3.5pt/:(1,-.2)@^{>}*!/3.5pt/:(1,+.2)@_{>}*!/7pt/:(1,-.2)@^{>}*!/7pt/:(1,+.2)@_{>}}
\newdir{ >>}{{}*!/8pt/@{|}*!/3.5pt/:(1,-.2)@^{>}*!/3.5pt/:(1,+.2)@_{>}}
\newdir{ |>}{{}*!/-3.5pt/@{|}*!/-8pt/:(1,-.2)@^{>}*!/-8pt/:(1,+.2)@_{>}}
\newdir{ >}{{}*!/-8pt/@{>}}
\newdir{>}{{}*:(1,-.2)@^{>}*:(1,+.2)@_{>}}
\newdir{<}{{}*:(1,+.2)@^{<}*:(1,-.2)@_{<}}

\newcommand{\n}{{\frak n}}
\newcommand{\g}{\frak g}
\newcommand{\q}{\frak q}

\newcommand{\Ker}{\ensuremath{\mathrm{Ker}}}
\newcommand{\Coker}{\ensuremath{\mathrm{Coker}}}

   \newcommand{\eq}{\frak q} 
  \newcommand{\eh}{\frak h}

\begin{document}

\centerline{\bf  ON $\Lie$-CENTRAL EXTENSIONS OF LEIBNIZ ALGEBRAS}

\bigskip
\bigskip
\centerline{\bf J. M. Casas$^{(1)}$ and E. Khmaladze$^{(2)}$}

\bigskip
\centerline{$^{(1)}$Dpto. Matemática Aplicada, Universidad de Vigo,  E. E. Forestal}
\centerline{Campus Universitario A Xunqueira, 36005 Pontevedra, Spain}
\centerline{ {E-mail address}: jmcasas@uvigo.es}
\medskip

\centerline{$^{(2)}$ A. Razmadze Mathematical Institute, Tbilisi State University}
\centerline{Tamarashvili Str. 6, 0177 Tbilisi, Georgia}
\centerline{E-mail address: {e.khmal@gmail.com}}

\date{}

\bigskip \bigskip \bigskip

{\bf Abstract:} Basing ourselves on the categorical notions of central extensions and commutators in the framework of semi-abelian categories relative to a Birkhoff subcategory, we study central extensions of Leibniz algebras with respect to the Birkhoff subcategory of Lie algebras, called $\Lie$-central extensions.
We obtain a six-term exact homology  sequence associated to a $\Lie$-central extension. This sequence, together with the relative commutators, allows us to characterize several classes of $\Lie$-central extensions, such as $\Lie$-trivial, $\Lie$-stem and $\Lie$-stem cover, to introduce and characterize $\Lie$-unicentral, $\Lie$-capable, $\Lie$-solvable and $\Lie$-nilpotent Leibniz algebras.
\bigskip

\bigskip


%---------------------------------------------------------------------------------------

\begin{section}{Introduction}

A general theory of central extensions relative to a chosen subcategory of a base category was introduced in \cite{JK}. Recently, in \cite{CVDL} were analyzed in details the categorical concepts of  central extensions, perfect objects and commutators in a semi-abelian category \cite{BB}, relative to a Birkhoff subcategory. Some examples (e.g. groups vs. abelian groups, Lie (Leibniz) algebras vs. vector spaces) are absolute, meaning that they are fitted into the relative case with respect to the subcategory of all abelian objects. In the absolute case, some results were already investigated in \cite{GTVDL}.

In this paper we deal with the 'non absolute' case: Leibniz algebras vs. Lie algebras. In particular, the goal of the present paper is to consider the relative notions of central extension and commutator when the semi-abelian category is {\Leib}, the category of Leibniz algebras, and its Birkhoff subcategory is $\Lie$, the category of Lie algebras, together with the Liezation functor $(-)_{\Lie} : {\Leib} \to {\Lie}$, which is left adjoint to the inclusion functor ${\Lie} \hookrightarrow {\Leib}$.

Under these circumstances, the notion of central extension relative to $\Lie$, so called $\Lie$-central extension, and the notion of commutator relative to $\Lie$, provide the necessary ingredients to introduce the notions of unicentrality, capability, solvability and nilpotency of Leibniz algebras relative to $\Lie$, all of them named with the prefix $\Lie$-.

Homological machinery relative to $\Lie$, coming from the categorical semi-abelian framework \cite{Ev, EVDL, EVDL1}, allows us  to characterize these new notions by means of the six-term exact sequence
\[
\Lien \otimes \Lieg_{\Lie} \longrightarrow HL^{\Lie}_2(\Lieg) \longrightarrow HL^{\Lie}_2(\Lieq)  \longrightarrow \Lien  \longrightarrow \Lieg_{\Lie} \longrightarrow \Lieq_{\Lie} \longrightarrow 0
\]
 associated to the $\Lie$-central extension $0 \to \Lien \to \Lieg  \to \Lieq \to 0$, where
  $HL^{\Lie}_{2}(-)$ denotes (the Hopf formula for) the second homology of a Leibniz algebra with coefficients in the Liezation functor $(-)_{\Lie}$  (see Subsections \ref{preliminaries} and \ref{six-term exact sequence}).

 In particular, we organize the paper as follows: in Section \ref{CB} we recall the necessary categorical background from \cite{CVDL, CVDL1}. Then we particularize in Subsection \ref{preliminaries} these notions for the category {\Leib} and its Birkhoff subcategory $\Lie$. So we recall the concepts of $\Lie$-perfect Leibniz algebras, $\Lie$-trivial and $\Lie$-central extensions of Leibniz algebras. In Subsection \ref{six-term exact sequence} we give the five-term exact sequence in homology with coefficients in the Liezation functor, associated to an extension of Leibniz algebras, and we obtain the above mentioned  six-term exact sequence associated to a $\Lie$-central extension.
In Subsection \ref{RCE} we introduce and study the notions of $\Lie$-stem extension and $\Lie$-stem cover of Leibniz algebras.
In Section \ref{Section-precise} we introduce the notions of $\Lie$-unicentral, $\Lie$-capable and precise $\Lie$-center of Leibniz algebras and investigate interrelationships between them.
%$C: Z_{\Lie}(\Lieq) \otimes {\Lieq}_{\Lie} \to HLR_2(\Lieq)$ in the six-term exact sequence associated to the $\Lie$-central extension $0 \to Z_{\Lie}(\Lieq) \%to \Lieq \to \Lieq / Z_{\Lie}(\Lieq) \to 0$. Here $Z_{\Lie}(\Lieq) = \{q \in \Lieq \mid [q,q']+[q',q]=0,\ {\rm for\ all}\ q' \in \Lieq \}$.
In Section \ref{nil}, by using the relative commutators, we introduce $\Lie$-central and  $\Lie$-abelian series of Leibniz algebras. Then we introduce and study the concepts of $\Lie$-solvability and $\Lie$-nilpotency of Leibniz algebras.
In Section \ref{HC} we obtain a homological characterization of $\Lie$-nilpotency of a Leibniz algebra, which is similar to the  Stallings' theorem for Leibniz algebras in the absolute case \cite{Ca}.

\end{section}

\begin{section}{Categorical background}\label{CB}
In this section we give an overview of needed definitions and results from \cite{CVDL, CVDL2, CVDL1}  on central extensions, perfect objects and commutators in   semi-abelian categories.
% \cite{BB}.

In their article \cite{JK}, Janelidze and Kelly introduced a general theory of  central extensions relative to a chosen subcategory $\B$ of the base category $\A$. This theory holds when $\B$ is a Birkhoff subcategory of a semi-abelian category $\A$ \cite{CVDL1}.

Recall from  \cite{JMT} that a category $\A$ is \emph{semi-abelian}  if it is pointed, Barr exact and Bourn protomodular with binary coproducts. A subcategory $\B$ of a semi-abelian category $\A$ is \emph{a Birkhoff subcategory}  if it is full, reflective and closed under subobjects and regular quotients.

From now on, we consider $\B$  a fixed Birkhoff subcategory of a semi-abelian category $\A$.
We denote by  $I :  {\A} \to {\B}$ the left adjoint to the inclusion functor $\B \hookrightarrow \A$, and write the components of its unit by $\eta_A : A \to I(A)$.

Then \emph{an extension} in $\A$ is a regular epimorphism in $\A$.

An extension $f\colon{B\twoheadrightarrow A}$ in $\A$ is said to be \emph{trivial (with respect to $\B$)} or \emph{$\B$-trivial} \cite{JK} when the induced square
 \begin{equation}\label{Birkhoff-Square}
\vcenter{\xymatrix{B \ar[r]^-{f} \ar[d]_-{\eta_{B}} & A \ar[d]^-{\eta_{A}}\\
I(B) \ar[r]_-{I(f)} & I(A)}}
\end{equation}
 is a pullback.

 An extension $f\colon{B\twoheadrightarrow A}$ in $\A$ is called \emph{central (with respect to $\B$)} or \emph{$\B$-central} \cite{JK} when either one of the projections $f_{0}$, $f_{1}$ in the kernel pair $(R[f],f_{0},f_{1})$ of $f$ is $\B$-trivial. That is to say, $f$ is $\B$-central if and only if in the diagram
\[
\xymatrix{R[f] \ar@<.5ex>[r]^-{f_{0}} \ar@<-.5ex>[r]_-{f_{1}} \ar[d]_-{\eta_{R[f]}} & B \ar[r]^-{f} \ar[d]^-{\eta_{B}} & A\\
IR[f] \ar@<.5ex>[r]^-{If_{0}} \ar@<-.5ex>[r]_-{If_{1}} & IB}
\]
either one of the left hand side squares is a pullback.

An object $P$ of $\A$ is called \emph{perfect with respect to $\B$} or \emph{$\B$-perfect} when $I(P)$ is the zero object $0$ of $\B$.

Since, in a semi-abelian category, a regular epimorphism is always the cokernel of its kernel, an appropriate notion of short exact sequence exists. Such will be any sequence ${K \overset{k}{\longrightarrow}  B \overset{f}{\longrightarrow} A}$ that satisfies $k=\Ker(f)$ and $f= {\Coker}(k)$. We denote this situation by $0\longrightarrow {K \overset{k}{\longrightarrow}  B \overset{f}{\longrightarrow} A}\longrightarrow 0$.

For any object $A$ of $\A$, the adjunction $\xymatrix@1{{\A} \ar@<1ex>[r]^-{I} \ar@{}[r]|-{\perp} & {\B} \ar@<1ex>[l]^-{\supset}}$ induces a short exact sequence
\[
\xymatrix{0 \ar[r] & [A, A]_{\B} \ar[r]^-{\mu_{A}} & A \ar[r]^-{\eta_{A}} & IA \ar[r] & 0.}
\]
Here the object $[A,A]_{\B}$, defined as the kernel of $\eta_{A}$, acts as \emph{a zero-dimensional commutator relative to $\B$}. Of course, $IA=A/[A,A]_{\B}$, so that $A$ is an object of~$\B$ if and only if $[A,A]_{\B}$ is zero. On the other hand, an object~$A$ of~$\A$ is called \emph{$\B$-perfect} when $[A,A]_{\B}=A$.

Let us remark that an extension $f\colon{B\twoheadrightarrow A}$ in $\A$ is $\B$-central if and only if the restrictions $[f_{0},f_{0}]_{\B}$, $[f_{1},f_{1}]_{\B}\colon{[R[f],R[f]]_{\B}\to [B,B]_{\B}}$ of the kernel pair projections $f_{0}$, $f_{1}\colon{R[f]\to B}$ coincide. This is the case precisely when $[f_{0},f_{0}]_{\B}$ and $[f_{1},f_{1}]_{\B}$ are isomorphisms, or, equivalently, when $\Ker[f_{0},f_{0}]_{\B}\colon{L_{1}[f]\to [R[f],R[f]]_{\B}}$  is zero. In the following diagram
\[
\xymatrix{& L_{1}[f] \ar[d]_-{\Ker [f_{0},f_{0}]_{\B}} \ar@{.>}[ddr] \\
0 \ar[r] & [R[f],R[f]]_{\B} \ar[r]^-{\mu_{R[f]}} \ar@<-.5ex>[d]_-{[f_{0},f_{0}]_{\B}} \ar@<.5ex>[d]^-{[f_{1},f_{1}]_{\B}} & R[f] \ar@<-.5ex>[d]_-{f_{0}} \ar@<.5ex>[d]^-{f_{1}} \ar[r]^-{\eta_{R[f]}} & I R[f] \ar@<-.5ex>[d]_-{I f_{0}} \ar@<.5ex>[d]^-{I f_{1}} \ar[r] & 0\\
0 \ar[r] & [B,B]_{\B} \ar[r]_-{\mu_{B}} & B \ar[r]_-{\eta_{B}} & IB \ar[r] & 0\; ,}
\]
through the composite $f_{1}\comp\mu_{R[f]}\comp\Ker[f_{0}]_{\B}$ the object $L_{1}[f]$ may be considered as a normal subobject of $B$. It acts as \emph{a one-dimensional commutator relative to $\B$} and, if $K$ denotes the kernel of $f$, it is usually written $[K,B]_{\B}$.

Let $A$ be an object of a semi-abelian category $\A$ with enough projectives  and $f\colon{B \twoheadrightarrow A}$ be a projective presentation with kernel~$K$. The induced objects
\[
\frac{[B,B]_{\B}}{[K,B]_{\B}}\qquad\qquad \text{and} \qquad\qquad \frac{K \cap [B,B]_{\B}}{[K,B]_{\B}}
\]
are independent of the chosen projective presentation of $A$ as explained for instance in~\cite{EVDL}. The object ${(K\cap [B,B]_{\B})}/{[K,B]_{\B}}$ is called (the Hopf formula for) the second homology object of $A$ (with coefficients in $\B$) and it is denoted by $H_{2}(A,I)$. We also write $H_{1}(A,I)$ for $I(A)$.

When $\A$ is a semi-abelian monadic category, the objects $H_{1}(A,I)$ and $H_{2}(A,I)$ may be computed using comonadic homology as in \cite{EVDL1} and they are fitted into the semi-abelian homology theory (see \cite{Ev}). Moreover, \cite[Theorem 5.9]{EVDL} states that, given a short exact sequence $0\longrightarrow {K {\longrightarrow}  B {\longrightarrow} A}\longrightarrow 0$,  there
is the following five-term exact sequence
\begin{equation}\label{five-term semi-ab}
\xymatrix{H_{2}(B,I) \ar[r] & H_{2}(A,I) \ar[r]^{\theta^{\ast}(B)} & \frac{K}{[K,B]_{\B}} \ar[r] & H_{1}(B,I) \ar[r] & H_{1}(A,I) \ar[r] & 0 \;.}
\end{equation}
\end{section}

%---------------------------------------------------------------------------------------

\begin{section}{Central extensions of Leibniz algebras with respect to Lie algebras} \label{central extensions}
 In this section we consider the particular case where $\A$ is the semi-abelian category {\Leib}\ of Leibniz algebras, and the Birkhoff subcategory $\B$ is $\Lie$, the category of Lie algebras.

%---------------------------------------------------------------------------------------

\begin{subsection}{Preliminary results on Leibniz algebras} \label{preliminaries}
We fix $\mathbb{K}$ as a ground field such that $\frac{1}{2} \in \mathbb{K}$. All vector spaces and tensor products are considered over $\mathbb{K}$.

A \emph{Leibniz algebra} \cite{Lo 1, Lo 2, LP} is a vector space $\eq$  equipped with a bilinear map $[-,-] : \Lieq \otimes \Lieq \to \Lieq$, usually called the \emph{Leibniz bracket} of $\eq$,  satisfying the \emph{Leibniz identity}:
\[
 [x,[y,z]]= [[x,y],z]-[[x,z],y], \quad x, y, z \in \Lieq.
\]
Leibniz algebras form a semi-abelian category \cite{CVDL, JMT}, denoted by {\Leib}, whose morphisms are linear maps that preserve the Leibniz bracket.

 A subalgebra ${\eh}$ of a Leibniz algebra ${\Lieq}$ is said to be \emph{left (resp. right) ideal} of ${\Lieq}$ if $ [h,q]\in {\eh}$  (resp.  $ [q,h]\in {\eh}$), for all $h \in {\eh}$, $q \in {\Lieq}$. If ${\eh}$ is both
left and right ideal, then ${\eh}$ is called \emph{two-sided ideal} of ${\Lieq}$. In this case $\Lieq/\Lieh$ naturally inherits a Leibniz algebra structure.

For a Leibniz algebra ${\Lieq}$, we denote by ${\Lieq}^{\rm ann}$ the subspace of ${\Lieq}$ spanned by all elements of the form $[x,x]$, $x \in \Lieq$. Further, we consider
\[
Z^r({\Lieq}) = \{a \in {\Lieq} \mid [x,a] = 0,\ x \in {\Lieq}\},\quad Z({\Lieq}) = \{a \in {\Lieq} \mid [x,a] = 0=[a,x],\ x \in {\Lieq}\}
\]
 and call the \emph{right center} and  \emph{center}  of $\Lieq$, respectively.  It is proved in \cite[Lemma 1.1]{KP} that both ${\Lieq}^{\rm ann}$ and $Z^r({\Lieq})$ are two-sided ideals of $\Lieq$. It is obvious that $Z({\Lieq})$ is also a  two-sided ideal of $\Lieq$.

Given a Leibniz algebra $\Lieq$, it is clear that the quotient ${\Lieq}_ {_{\rm Lie}}=\Lieq/{\Lieq}^{\rm ann}$ is a Lie algebra. This defines the so-called  \emph{Liezation functor} $(-)_{\Lie} : {\Leib} \to {\Lie}$, which assigns to a Leibniz algebra $\Lieq$ the Lie algebra ${\Lieq}_{_{\rm Lie}}$. Moreover,
the canonical epimorphism  ${\Lieq} \twoheadrightarrow {\Lieq}_ {_{\rm Lie}}$ is universal among all homomorphisms from $\Lieq$ to a Lie algebra, implying that the Liezation functor is left adjoint to the inclusion functor $ {\Lie} \hookrightarrow {\Leib} $.

 It is an easy task to check that the category $\Lie$ is a Birkhoff subcategory of $\Leib$.
Focusing our attention in the  adjoint pair
\begin{equation}
\xymatrix@1{{\Leib} \ar@<1ex>[r]^-{(-)_{\Lie}} \ar@{}[r]|-{\perp} & {\Lie} \ar@<1ex>[l]^-{\supset}},
\end{equation}
we  particularize the general theory described in Section \ref{CB} to the case when $\A$ is the category {\Leib}, $\B$ is its Birkhoff subcategory $\Lie$
and the functor $I$ is precisely the Liezation functor $(-)_{\Lie}$. First let us remark the following

\begin{Rem}
A Leibniz algebra $\q$ is $\Lie$-perfect if $\q_{\Lie}=0$, that is $\q \cong \q^{\rm ann}$. It follows by \cite[Lemma 3.1]{CILL} that any $\Lie$-perfect Leibniz algebra $\q$ is the trivial one. Then, by \cite[Theorem 3.5]{CVDL1}, we get that a Leibniz algebra $\q$ admits a universal $\Lie$-central extension if and only if $\q=0$.
\end{Rem}

It is clear that an extension in $\Leib$ is just an epimorphism $f\colon{\Lieg \twoheadrightarrow \Lieq}$ of Leibniz algebras, which is a $\Lie$-trivial extension if and only if it induces an isomorphism $\Lieg^{\rm ann} {\cong}\Lieq^{\rm ann}$.
 \begin{Ex}
 Let $\Lieg$ and $\Lieq$  be three and two-dimensional (as vector spaces) Leibniz algebras with $\K$-linear bases $\{a_1, a_2, a_3 \}$ and $\{e_1, e_2\}$, respectively, with the Leibniz brackets given respectively by  $[a_1,a_3]=a_1$ and $[e_1,e_2]=e_1$ and zero elsewhere (see \cite{CILL, Cu}).
 Consider the homomorphism of Leibniz algebras $f : \Lieg \to \Lieq$ defined by $f(a_1)=e_1$, $f(a_2)=0$, $f(a_3)=e_2$. Obviously $f$ is surjective, $\Lieg^{\rm ann}=\langle \{a_1 \} \rangle$ and $\Lieq^{\rm ann}=\langle \{e_1 \} \rangle$, so they are isomorphic through $f$. Consequently $f : \Lieg \twoheadrightarrow \Lieq$ is a $\Lie$-trivial extension.
\end{Ex}

For a Leibniz algebra {\Lieq} and two-sided ideals  ${\Liem}$ and ${\Lien}$ of  ${\Lieq}$, we put
\[
C_{\Lieq}^{\Lie}({\Liem} , {\Lien}) = \{q \in {\Lieq} \mid  \; [q, m] + [m,q] \in {\Lien}, \; \text{for all} \;
m \in {\Liem} \} \; .
\]

\noindent Further, we denote by $[\Liem,\Lien]_{\Lie}$ the subspace  of $\Lieq$ spanned by all elements of the form $[m,n]+[n,m]$, $m \in \Liem$, $n \in \Lien$.

\begin{Le}\label{Lemma3.3} Let $\Lieq$ be a Leibniz algebra and $\Liem$, $\Lien$ be two-sided ideals of $\Lieq$. Then both $C_{\Lieq}^{\Lie}({\Liem} , {\Lien})$  and $[\Liem,\Lien]_{\Lie}$ are  two-sided ideals of $\Lieq$. Moreover, $ Z(\Lieq)\subseteq C_{\Lieq}^{\Lie}({\Liem} , {\Lien}) $ and $[\Liem,\Lien]_{\Lie}\subseteq Z^r(\Lieq)$.
\end{Le}
{\it Proof.}
It is clear that both  $C_{\Lieq}^{\Lie}({\Liem} , {\Lien})$ and $[\Liem,\Lien]_{\Lie}$ are subspaces of $\Lieq$.

Take any elements $q\in C_{\Lieq}^{\Lie}({\Liem} , {\Lien})$ and $x\in \Lieq$. Then, for any $m\in \Liem$, we have
\begin{align*}
[[q,x],m]+[m,[q,x]]&= -[q, [m,x]]+[[q,m],x]+[[m,q],x]-[[m,x],q] \\
&= -\big([q, [m,x]]+[[m,x],q]\big) + \big[[q,m]+[m,q] , x\big] \in \Lien
\end{align*}
and, if we denote $-\big([q, [m,x]]+[[m,x],q]\big) + \big[[q,m]+[m,q] , x\big]$ by $n$, we get
\begin{align*}
[[x,q],m]+[m,[x,q]]&= [[x,q],m]-[m,[q,x]]
 = [[x,q],m]- n + [[q,x],m]\\
&=\big[[x,q]+[q,x], m\big]-n \in \Lien.
\end{align*}
Consequently, $C_{\Lieq}^{\Lie}({\Liem} , {\Lien})$  is a two-sided ideal of $\Lieq$. The inclusion $Z(\Lieq) \subseteq C_{\Lieq}^{\Lie}({\Liem} , {\Lien})$ is obvious.

Now, for any $m\in \Liem$, $n\in \Lien$ and $q\in \Lieq$, we have
\[
[[m,n]+[n,m],q]=[m,[n,q]]+[[n,q],m]+[n,[m,q]]+[[m,q],n] \in [\Liem,\Lien]_{\Lie},
\]
and
\[
[q,[m,n]+[n,m]]=0.
\]
Consequently, $[\Liem,\Lien]_{\Lie}$ is a two-sided ideal of $\Lieq$ contained in $Z^r(\Lieq)$.
\rdg

In particular, the two-sided ideal $C_{\Lieq}^{\Lie}(\Lieq , 0)$ is the \emph{$\Lie$-center} of the Leibniz algebra $\Lieq$ and it will be denoted by $Z_{\Lie}(\Lieq)$, that is,

\[
Z_{\Lie}(\Lieq) =  \{ z\in \Lieq\,|\,\text{$[q,z]+[z,q]=0$ for all $q\in \Lieq$}\}.
\]

The following proposition is an immediate consequence of the discussion in \cite[Example 1.9]{CVDL}.
\begin{Pro}
Given an extension of Leibniz algebras $f : \Lieg \twoheadrightarrow \Lieq$ with ${\Lien} =  \Ker(f)$, the following conditions are equivalent:
\begin{enumerate}
\item[$(a)$] $f : \Lieg \twoheadrightarrow \Lieq$ is $\Lie$-central;
\item[$(b)$] $\Lien \subseteq Z_{\Lie}(\Lieq)$;
\item[$(c)$] $[\Lien,\Lieg]_{\Lie} = 0$.
\end{enumerate}
\end{Pro}

 \begin{Rem}
 Obviously every central extension $f : \Lieg \twoheadrightarrow \Lieq$  (i.e. $\Ker (f) \subseteq Z(\Lieg)$) is a $\Lie$-central extension, but the converse is not true in general as the following example shows: consider the three-dimensional Leibniz algebra $\Lieg$  with $\K$-linear basis $\{a_1,a_2,a_3\}$, the Leibniz bracket given  by $[a_1,a_3]=a_1$ and zero elsewhere, and the two-dimensional abelian Leibniz algebra $\Lieq$ with $\K$-linear basis $\{e_1, e_2\}$. Then the homomorphism of Leibniz algebras $f : \Lieg \to \Lieq$ given by $f(a_1)=0, f(a_2)=e_1, f(a_3)=e_2$ is surjective, ${\rm Ker}(f) = \langle \{a_1 \} \rangle$, $Z(\Lieg) =\langle \{a_2 \} \rangle$ and $Z_{\Lie}(\Lieg)=\langle \{a_1, a_2 \} \rangle$. Consequently the extension $f : \Lieg \twoheadrightarrow \Lieq$ is not central, but it is $\Lie$-central.
 \end{Rem}

\end{subsection}

%---------------------------------------------------------------------------------------

\begin{subsection}{Six-term exact sequence}\label{six-term exact sequence}

In what follows, given a Leibniz algebra $\Lieg$, we write $HL^{\Lie}_2(\Lieg)$ for $H_2(\Lieg, (-)_{\Lie})$ and  $HL^{\Lie}_1(\Lieg)$)  for $H_1(\Lieg, (-)_{\Lie})$. Thus, the five-term exact sequence (\ref{five-term semi-ab}) associated to an  extension of Leibniz algebras $f\colon{\Lieg \twoheadrightarrow \Lieq}$, with $\Lien = {\rm Ker}(f)$,  turns to
\begin{equation}\label{five-term trivial} HL^{\Lie}_2(\Lieg) \longrightarrow HL^{\Lie}_2(\Lieq) \stackrel{\theta^{\ast}(\Lieg)} \longrightarrow \frac{\Lien}{[\Lien,\Lieg]_{\Lie}} \longrightarrow HL^{\Lie}_1(\Lieg) \longrightarrow HL^{\Lie}_1(\Lieq) \longrightarrow 0.
\end{equation}
Furthermore, we have the following proposition.
\begin{Pro}\label{Proposition3.9}
Given a $\Lie$-central extension $f:\Lieg\twoheadrightarrow \Lieq$ of Leibniz algebras with $\Lien=\Ker(f)$, there is a six-term exact sequence
\begin{equation}\label{six-term}
 \Lien\otimes \Lieg_{\Lie} \longrightarrow HL^{\Lie}_2(\Lieg) \longrightarrow HL^{\Lie}_2(\q) \stackrel{\theta^{\ast}(\Lieg)}\longrightarrow \frak{n} \longrightarrow HL^{\Lie}_1(\Lieg) \longrightarrow HL^{\Lie}_1(\q) \longrightarrow 0.
\end{equation}
\end{Pro}
{\it{Proof.}} Any free presentation of $\Lieg$, $0 \to \Lier \to \Lief \stackrel{\rho} \to \Lieg \to 0$, produces a free presentation of $\Lieq$,  $0 \to \Lies \to \Lief \stackrel{\tau=f \circ \rho} \to \Lieq \to 0$, and the following commutative diagram of Leibniz algebras
\begin{equation} \label{free presentations}
\xymatrix{& & 0 \ar[d] & 0 \ar[ld]  &\\
&  & \Lier \ar[d]\ar[ld] & &\\
0 \ar[r] & \frak{s}\ \ar[r] \ar[d] & \Lief \ar[d]^{\rho} \ar[rd]^{\tau=f \circ \rho} & &\\
0 \ar[r]& \Lien \ar[r] \ar[d]& \Lieg \ar[r]^{f} \ar[d]& \Lieq \ar[r] \ar[dr] & 0 \\
 & 0 & 0 &  & 0
}  \end{equation}
Thus we have
\[
\Ker\left(HL^{\Lie}_2(\Lieg) \to HL^{\Lie}_2(\q)\right) \approx \Ker\left( \frac{\Lier \cap [\Lief,\Lief]_{\Lie}}{[\Lier,\Lief]_{\Lie}} \to \frac{\Lies \cap [\Lief,\Lief]_{\Lie}}{[\Lies,\Lief]_{\Lie}}\right)\approx \frac{[\Lies,\Lief]_{\Lie}}{[\Lier,\Lief]_{\Lie}}.
\]

When $f$ is a {\Lie}-central extension, $[\Lien,\Lieg]_{\Lie}=0$, and  by the exact sequence (\ref{five-term trivial}) we get the following exact sequence
\begin{equation}\label{five-term trivial 1}
0\longrightarrow \frac{[\Lies,\Lief]_{\Lie}}{[\Lier,\Lief]_{\Lie}}\longrightarrow \frac{\frak{r} \cap [\frak{f}, \frak{f}]_{\Lie}}{[\frak{r},\frak{f}]_{\Lie}} \longrightarrow \frac{\frak{s} \cap [\frak{f}, \frak{f}]_{\Lie}}{[\frak{s},\frak{f}]_{\Lie}} \stackrel{\theta^{\ast}(\Lieg)} \longrightarrow \Lien \longrightarrow \Lieg_{\Lie} \longrightarrow \Lieq_{\Lie} \longrightarrow 0.
\end{equation}
%the sequence (\ref{five-term trivial 1}) can be rewritten  as follows:
%\begin{equation}\label{five-term} HL^{\Lie}_2(\Lieg) \longrightarrow HL^{\Lie}_2(\q) \stackrel{\theta^{\ast}(\Lieg)} \longrightarrow \frak{n} \longrightarrow %HL^{\Lie}_1(\Lieg) \longrightarrow HL^{\Lie}_1(\q) \longrightarrow 0.
%\end{equation}

Now consider $\sigma: \Lien \otimes \Lieg_{\Lie} \to \frac{[\Lies,\Lief]_{\Lie}}{[\Lier,\Lief]_{\Lie}}$ given by $\sigma(n \otimes \overline{x}) = [s,f]+[f,s]+ [\Lier,\Lief]_{\Lie}$, where $s\in \Lier$ and $f\in \Lief$  such that $\rho(s)=n$ and $\rho(f)=x$. Again using the condition $[\Lien,\Lieg]_{\Lie}=0$, it is easy to check that  $\sigma$ is a well-defined epimorphism. Then the assertion follows. \rdg

\end{subsection}

%---------------------------------------------------------------------------------------

%---------------------------------------------------------------------------------------

\begin{subsection}{Some properties of $\Lie$-central extensions} \label{RCE}

In this subsection we establish some properties of $\Lie$-central extensions of Leibniz algebras by using
the six-term exact sequence (\ref{six-term}) associated to the given extension.

\begin{Pro} \label{Lie trivial}
Let $f\colon{\Lieg \twoheadrightarrow \Lieq}$ be an extension of Leibniz algebras with $\Lien = \Ker (f)$. The
following statements are equivalent:
\begin{enumerate}
\item[(a)] $f\colon{\Lieg \twoheadrightarrow \Lieq}$ is a $\Lie$-trivial extension.
\item[(b)] $\Lien \cap \Lieg^{\rm ann} =0$.
 \item[(c)] $0 \to {\Lien}  \to {\Lieg}_{_{\Lie}} \to
{\Lieq}_{_{\Lie}} \to 0$ is exact in $\Lie$.
\end{enumerate}
\end{Pro}
{\it Proof.} The equivalences follow from the definition of $\Lie$-trivial extension  and by the $3 \times 3$ Lemma \cite{B2} applied to the following diagram:
\[ \xymatrix{
 & 0 \ar[d] & 0 \ar[d] & 0 \ar[d] & \\
0 \ar[r] & \n \cap \g^{\rm ann} \ar[r] \ar[d] & \g^{\rm ann} \ar[r] \ar[d] & \q^{\rm ann} \ar[r] \ar[d] & 0\\
0 \ar[r] & \n \ar[r] \ar[d] & \g \ar[r] \ar[d] & \q \ar[r] \ar[d] & 0\\
0 \ar[r] & \frac{\n}{\n \cap \g^{\rm ann}} \ar[r] \ar[d] & HL^{\Lie}_1(\g) \ar[r] \ar[d] & HL^{\Lie}_1(\q) \ar[r] \ar[d] & 0\\
 & 0  & 0 & 0  &
} \]
 \rdg

\begin{Pro}
Let $f\colon{\Lieg \twoheadrightarrow \Lieq}$ be a $\Lie$-trivial extension with $\Lien = \Ker (f)$, then $\theta^{\ast}(\Lieg) : HL^{\Lie}_2(\Lieq) \to \Lien$ is the zero map and $\Lien \otimes \Lieg_{\Lie} \to HL^{\Lie}_2(\Lieg) \to HL^{\Lie}_2(\Lieq) \to 0$ is exact.
\end{Pro}
{\it Proof.} Since $[\n,\g]_{\Lie} \subseteq \n \cap \g^{\rm ann}$, by Proposition \ref{Lie trivial} we have that $[\n,\g]_{\Lie} = 0$, then $\Im (\theta^{\ast}(\g))=0$ in sequence (\ref{six-term}). \rdg

\begin{Co}
Let $f\colon{\Lieg \twoheadrightarrow \Lieq}$  be a $\Lie$-trivial extension with \Lien = \Ker($f$), then $\Lien = Z_{\Lie}(\Lien)$.
\end{Co}
{\it Proof.} Since $\n \cap \g^{\rm ann}=0$, we have $[n,n']+[n',n]=[n+n',n+n']=0$ for any $n, n' \in \n$, which means that $\n \subseteq Z_{\rm Lie}(\n)$. \rdg
\bigskip

%Now we consider a  $\Lie$-central extension $f : \Lieg \twoheadrightarrow\Lieq$ with \Lien = \Ker($f$).

\begin{De}
A  $\Lie$-central extension $f : \Lieg \twoheadrightarrow \Lieq$ is said to be:
\begin{enumerate}
\item[(a)] a $\Lie$-stem extension if $\Lieg_{\Lie} \cong \Lieq_{\Lie}$.
\item[(b)] a $\Lie$-stem cover if $\Lieg_{\Lie} \cong \Lieq_{\Lie}$ and the induced map $HL^{\Lie}_2(\Lieg) \to HL^{\Lie}_2(\Lieq)$ is the zero map.
\end{enumerate}
\end{De}

\begin{Ex}\ \label{Example3.14}
\begin{enumerate}
\item[(a)]  Let $\Lieg$ and $\Lieq$ be  three and two-dimensional Leibniz algebra with $\K$-linear bases $\{a_1, a_2, a_3\}$  and $\{e_1,e_2 \}$, with the Leibniz brackets given respectively by $[a_2,a_3] = - [a_3,a_2] =a_2$, $[a_3,a_3]= a_1$ and $[e_1,e_2] = - [e_2, e_1] = -e_1$ and zero elsewhere. Then the surjective homomorphism of Leibniz algebras $f : \Lieg \to \Lieq$ defined by $f(a_1)=0$, $f(a_2)=e_1$ and $f(a_3)=e_2$ is a $\Lie$-stem extension.
\item[(b)] Let $0 \to \Lies \to \Lief \stackrel{\tau} \to \Lieq \to 0$ be a free presentation of a Leibniz algebra $\Lieq$. Then $[\Lies, \Lief]_{\Lie}$ is a two-sided ideal of $\Lief$ and the $3 \times 3$ Lemma \cite{B2}  provides the sequence $0 \to \Lies/[\Lies, \Lief]_{\Lie} \to \Lief/[\Lies, \Lief]_{\Lie} \to \Lieq \to 0$ which is a $\Lie$-stem cover of $\Lieq$.
\end{enumerate}
\end{Ex}

\begin{Pro}
For a $\Lie$-central extension $f : \Lieg \twoheadrightarrow \Lieq$, with $\Lien = {\rm Ker}(f)$, the following statements are equivalent:
\begin{enumerate}
\item[(a)] $f : \Lieg \twoheadrightarrow \Lieq$ is a $\Lie$-stem extension.
\item[(b)] The induced map $\Lien \to HL^{\Lie}_1(\Lieq)$ is the zero map.
\item[(c)] $\theta^{\ast}(\Lieg) : HL^{\Lie}_2(\Lieq) \to \Lien$ is an epimorphism.
\item[(d)] The following sequence $\n \otimes \Lieg_{\Lie} \to HL^{\Lie}_2(\Lieg) \to HL^{\Lie}_2(\q) \stackrel{\theta^{\ast}(\Lieg)} \to \frak{n} \to 0$ is exact.
\item[(e)]$\Lien \subseteq \Lieg^{\rm ann}$.
\end{enumerate}
\end{Pro}
{\it Proof.} The equivalences between (a), (b), (c) and (d) follow from the exact sequence (\ref{six-term}).
The equivalence between (a) and (e) is a consequence of the following $3 \times 3$ diagram:
\[ \xymatrix{
 & 0 \ar[d] & 0 \ar[d] & 0 \ar[d] & \\
0 \ar[r] & \n \cap \g^{\rm ann} \ar[r] \ar[d] & \g^{\rm ann} \ar[r] \ar[d] & \q^{\rm ann} \ar[r] \ar[d] & 0\\
0 \ar[r] & \n \ar[r] \ar[d] & \g \ar[r] \ar[d] & \q \ar[r] \ar[d] & 0\\
0 \ar[r] & 0 \ar[r] \ar[d] & HL^{\Lie}_1(\g) \ar@{=}[r] \ar[d] & HL^{\Lie}_1(\q) \ar[r] \ar[d] & 0\\
 & 0  & 0 & 0  &
} \] \rdg

\begin{Pro} \label{stem cover}
For a  $\Lie$-central extension $f : \Lieg \twoheadrightarrow \Lieq$  the following statements are equivalent:
\begin{enumerate}
\item[(a)] $f : \Lieg \twoheadrightarrow \Lieq$  is a $\Lie$-stem cover.
\item[(b)] $\theta^{\ast}(\Lieg) : HL^{\Lie}_2(\Lieq) \to \n$ is an isomorphism.
\end{enumerate}
\end{Pro}
{\it Proof.} This is a direct consequence of the six-term exact  sequence (\ref{six-term}). \rdg

\end{subsection}

\end{section}

%---------------------------------------------------------------------------------------

\begin{section}{The precise $\Lie$-center of a Leibniz algebra} \label{Section-precise}
In this section we introduce the notions of $\Lie$-unicentrality, $\Lie$-capability and  precise $\Lie$-center of a Leibniz algebra and
analyze the relationships between them.
%% sequence (\ref{six-term}).

\begin{De}
A Leibniz algebra $\Lieq$ is said to be $\Lie$-unicentral if every $\Lie$-central extension $f : \Lieg \twoheadrightarrow \Lieq$ satisfies $f\left( Z_{\Lie}(\Lieg)\right) = Z_{\Lie}(\Lieq)$, that is, the following diagram with exact rows
\[\xymatrix{
 &{\Ker}(f) \ar@{^{(}->}[r] \ar@{=}[d] &Z_{\Lie}(\Lieg) \ar@{^{(}->}[d] \ar@{>>}[r] & Z_{\Lie}(\Lieq) \ar@{^{(}->}[d] \\
 &{\Ker}(f)\ \ar@{^{(}->}[r]\  & \Lieg \ar@{>>}[r]^{f} & \Lieq
}\]
is commutative.
\end{De}

 \begin{De}
 A Leibniz algebra $\Lieq$ is said to be $\Lie$-capable if there exists a $\Lie$-central extension
 \[
 0 \longrightarrow Z_{\Lie}(\Lieg) \longrightarrow \Lieg \stackrel{f} \longrightarrow \Lieq \longrightarrow 0.
 \]
 \end{De}

\begin{De}
The precise $\Lie$-center $Z_{\Lie}^{\ast}(\Lieq)$ of a Leibniz algebra $\Lieq$ is the intersection of all two-sided ideals $f(Z_{\Lie}(\Lieg))$, where $f : \Lieg \twoheadrightarrow \Lieq$ is a $\Lie$-central extension.
\end{De}

\begin{Rem}\
\begin{enumerate}
\item[(a)] $Z_{\Lie}^{\ast}(\Lieq) \subseteq Z_{\Lie}(\Lieq)$.

\item[(b)] $Z_{\Lie}^{\ast}(\Lieq) = Z_{\Lie}(\Lieq)$ if and only if $f(Z_{\Lie}(\Lieg)) = Z_{\Lie}(\Lieq)$ for every $\Lie$-central extension of Leibniz algebras $f : \Lieg \twoheadrightarrow \Lieq$, or equivalently, if and only if $\Lieq$ is $\Lie$-unicentral.
\end{enumerate}
\end{Rem}

%Now we  characterize $\Lie$-capability in terms of the cancelation  of the $\Lie$-precise center.
Given a free presentation $0 \to \Lies \to \Lief \stackrel{\tau} \to \Lieq \to 0$ of the Leibniz algebra $\Lieq$, consider the $\Lie$-central extension ($\Lie$-stem cover of $\Lieq$)
\[
0 \longrightarrow \frac{\Lies}{[\Lies,\Lief]_{\Lie}} \longrightarrow \frac{\Lief}{[\Lies,\Lief]_{\Lie}} \stackrel{\overline{\tau}} \longrightarrow \Lieq \longrightarrow 0
\]
as in Example \ref{Example3.14} (b). Then we have

\begin{Le} \label{precise}
$Z_{\Lie}^{\ast}(\Lieq) = \overline{\tau} \left( Z_{\Lie} \left( \frac{\Lief}{[\Lies, \Lief]_{\Lie}} \right) \right)$.
\end{Le}
{\it Proof.} We need to show that, for any $\Lie$-central extension $0 \to \Liea \to \Lieh \stackrel{\varphi}\to \Lieq \to 0$, the inclusion $\overline{\tau} \left( Z_{\Lie} \left( \frac{\Lief}{[\Lies, \Lief]_{\Lie}} \right) \right) \subseteq \varphi \left(Z_{\Lie}(\Lieh)\right)$ holds.

Since $\Lief$ is a free Leibniz algebra, there exists (uniquely defined epimorphism) $\alpha : \Lief \to \Lieh$ such that $\varphi  \circ \alpha = \tau$. Then $\alpha(\Lies) \subseteq \Liea$ and $\alpha\left( [\Lies,\Lief]_{\Lie} \right) \subseteq [\Liea, \Lieh]_{\rm Lie} =0$. Hence, $\alpha$ induces $\overline\alpha: \frac{\Lief}{[\Lies,\Lief]_{\Lie}} \to \Lieh$ such that $\overline\alpha \circ \pi = \alpha$, where $\pi:\Lief \twoheadrightarrow \frac{\Lief}{[\Lies,\Lief]_{\Lie}}$ is the natural projection. It is straightforward to see that
$\overline\alpha \left( Z_{\Lie} \left( \frac{\Lief}{[\Lies,\Lief]_{\Lie}} \right) \right) \subseteq Z_{\Lie}(\Lieh)$. Now, since $\overline{\tau}\circ \pi = \tau = \varphi \circ \alpha = \varphi \circ \overline\alpha \circ \pi$, it follows that
\[\overline{\tau} \left( Z_{\Lie} \left( \frac{\Lief}{[\Lies,\Lief]_{\Lie}} \right) \right) = (\varphi \circ \overline\alpha) \left( Z_{\Lie} \left( \frac{\Lief}{[\Lies,\Lief]_{\Lie}} \right) \right) \subseteq \varphi \left( Z_{\Lie} (\Lieh) \right).
\]
 \rdg

\begin{Co}\label{Corollary4.6}
$Z_{\Lie}^{\ast}(\Lieq) = 0$ if and only if $\Lieq$ is a $\Lie$-capable Leibniz algebra.
\end{Co}
{\it Proof.} If $\Lieq$ is a $\Lie$-capable Leibniz algebra, then there exists a $\Lie$-central extension $0 \to Z_{\Lie}(\Lieg) \to \Lieg \stackrel{f} \to \Lieq \to 0$, then $Z_{\Lie}^{\ast}(\Lieq) \subseteq f\left( Z_{\Lie}(\Lieg) \right) = 0$.

Conversely, if $Z_{\Lie}^{\ast}(\Lieq)=0$, for any free presentation $0 \to  \Lies \to \Lief \stackrel{\tau} \to \Lieq \to 0$ of $\Lieq$, we have
$ \overline{\tau} \left( Z_{\Lie} \left( \frac{\Lief}{[\Lies,\Lief]_{\Lie}} \right) \right)=0$  by Lemma \ref{precise}. Then
\[
0 \longrightarrow Z_{\Lie} \left( \frac{\Lief}{[\Lies,\Lief]_{\Lie}} \right) \longrightarrow \frac{\Lief}{[\Lies,\Lief]_{\Lie}} \stackrel{\overline{\tau}}\longrightarrow \Lieq \longrightarrow 0
\]
 is a $\Lie$-central extension. \rdg

\begin{Le} \label{epim}
Let $\pi : \Lieg \twoheadrightarrow \Lieq$ be an epimorphism of Leibniz algebras, then $\pi \left( Z_{\Lie}^{\ast}(\Lieg) \right) \subseteq Z_{\Lie}^{\ast}(\Lieq)$.
\end{Le}
{\it Proof.}  For any $\Lie$-central extension $0 \to \Lien \to \Lieh \stackrel{\varphi} \to \Lieq \to 0$ of $\Lieq$, consider the pull-back diagram over $\pi$ \[ \xymatrix{
0 \ar[r] & \Lien \ar[r] \ar@{=}[d] &\Lieg \times_{\Lieq} \Lieh \ar[r]^{\quad \overline{\varphi}} \ar@{>>}[d]^{\overline{\pi}} & \Lieg \ar[r] \ar@{>>}[d]^{\pi}& 0\\
0 \ar[r] & \Lien \ar[r] & \Lieh \ar[r]^{\varphi}  & \Lieq \ar[r]& 0 \ .
}\]
Clearly the upper row is again a $\Lie$-central extension of $\Lieq$. Then we have $\pi \left( Z_{\Lie}^{\ast}(\Lieg) \right) \subseteq \pi \circ \overline{\varphi} \left( Z_{\Lie}( \Lieg \times_{\Lieq} \Lieh) \right) = \varphi \circ \overline{\pi} \left( Z_{\Lie}( \Lieg \times_{\Lieq} \Lieh ) \right) \subseteq \varphi \left( Z_{\Lie}(\Lieh) \right)$. This implies that  $\pi \left( Z_{\Lie}^{\ast}(\Lieg) \right) \subseteq \underset{\varphi}{\bigcap} \varphi \left( Z_{\Lie}(\Lieh) \right) = Z_{\Lie}^{\ast}(\Lieq)$. \rdg

\begin{Pro}\
\begin{enumerate}
\item[(a)] $Z_{\Lie}^{\ast}(\Lieq)$ is the smallest two-sided ideal $\Lien$ of the Leibniz algebra $\Lieq$ such that $\Lieq/\Lien$ is $\Lie$-capable. In particular, $\Lieq/Z_{\Lie}^{\ast}(\Lieq)$ is $\Lie$-capable.
    \item[(b)] Let $\Lien$ be a two-sided ideal of a Leibniz algebra $\Lieq$ such that $\Lien \cap Z_{\Lie}^{\ast}(\Lieq)=0$. If $\Lieq/\Lien$ is $\Lie$-capable, then $\Lieq$ is $\Lie$-capable as well.
\end{enumerate}
\end{Pro}
{\it Proof.} (a) Given a two-sided ideal $\Lien$ of $\Lieq$ such that $\Lieq/\Lien$ is $\Lie$-capable, consider the epimorphism $\pi: \Lieq  \twoheadrightarrow \Lieq/\Lien $. By Lemma \ref{epim} and Corollary \ref{Corollary4.6} we have $\pi \left( Z_{\Lie}^{\ast}(\Lieq) \right) \subseteq Z_{\Lie}^{\ast}(\Lieq /\Lien) =0$, then $Z_{\Lie}^{\ast}(\Lieq) \subseteq {\rm Ker}(\pi) = \Lien$.

(b) By (a), $Z_{\Lie}^{\ast}(\Lieq) \subseteq  \Lien$. Since $Z_{\Lie}^{\ast}(\Lieq) \cap \Lien = 0$, then $Z_{\Lie}^{\ast}(\Lieq) = 0$. Then Corollary \ref{Corollary4.6} completes the proof. \rdg

\begin{Th}
Let $\Liea$ be a two-sided ideal of a Leibniz algebra $\Lieq$ such that $\Liea \subseteq Z_{\Lie}(\Lieq)$. Then
$\Liea \subseteq Z_{\Lie}^{\ast}(\Lieq)$ if and only if the map $C : \Liea \otimes \Lieq_{\Lie} \to HL^{\Lie}_2(\Lieq)$ in sequence $\left( \ref{six-term} \right)$ associated to the $\Lie$-central extension $0 \to \Liea \to \Lieq \stackrel{\pi} \to \Lieq/\Liea \to 0$ is the zero map.
\end{Th}
{\it Proof.}  Consider the free presentations $0 \to \Lier \to \Lief \stackrel{\rho} \to \Lieq \to 0$ and $0 \to \Lies \to \Lief \stackrel{\pi \circ \rho} \to \Lieq/\Liea \to 0$. We know from the proof of Proposition \ref{Proposition3.9} that Im$(C) \cong \frac{[\Lies, \Lief]_{\Lie}}{[\Lier, \Lief]_{\Lie}}$, then by Proposition \ref{precise} we have the following commutative diagram:
\[ \xymatrix{
& & Z_{\Lie} \left( \frac{\Lief}{[\Lier,\Lief]_{\Lie}} \right) \ar[r] \ar@{^{(}->}[d]&  Z_{\Lie}^{\ast}(\Lieq) \ar@{^{(}->}[d]& \\
0 \ar[r] & \frac{\Lier}{[\Lier,\Lief]_{\Lie}} \ar[r] \ar[ur] & \frac{\Lief}{[\Lier,\Lief]_{\Lie}} \ar[r]^{\overline{\rho}} \ar@{>>}[d]^{\epsilon} & \Lieq \ar[r] \ar@{>>}[d]^{\gamma} & 0\\
 & & \bullet \ar[r]^{\sim \ \ \ } & \Lieq / Z_{\Lie}^{\ast}(\Lieq) &
}\]
Hence $C=0 \Leftrightarrow \frac{[\Lies, \Lief]_{\Lie}}{[\Lier, \Lief]_{\Lie}}=0 \Leftrightarrow   \frac{\Lies}{[\Lier, \Lief]_{\Lie}}             \subseteq Z_{\Lie} \left( \frac{ \Lief}{[\Lier, \Lief]_{\Lie}} \right) \Leftrightarrow \gamma  \circ \overline{\rho} \left( \frac{\Lies}{[\Lier, \Lief]_{\Lie}} \right) = \epsilon \left( \frac{\Lies}{[\Lier, \Lief]_{\Lie}} \right) =0 \Leftrightarrow \Liea = \rho(\Lies) = \overline{\rho} \left(\frac{\Lies}{[\Lier, \Lief]_{\Lie}} \right) \subseteq {\rm Ker}(\gamma) = Z_{\Lie}^{\ast}(\Lieq)$. \rdg

\begin{Co}\label{charact}
For a Leibniz algebra $\Lieq$ the following statements are equivalent:
\begin{enumerate}
\item[(a)] $\Lieq$ is $\Lie$-unicentral.
\item[(b)] The map $C : Z_{\Lie}(\Lieq) \otimes \Lieq_{\Lie} \to HL^{\Lie}_2(\Lieq)$ in the sequence $\left( \ref{six-term} \right)$ associated to the $\Lie$-central extension $0 \to Z_{\Lie}(\Lieq) \to \Lieq \to \Lieq/Z_{\Lie}(\Lieq) \to 0$ is the zero map.

    \item[(c)] The canonical homomorphism $HL^{\Lie}_2(\Lieq) \to HL^{\Lie}_2(\Lieq/Z_{\Lie}(\Lieq))$ is injective.
\end{enumerate}
\end{Co}
{\it Proof.} This is a consequence of the exactness of the sequence (\ref{six-term}) associated to the $\Lie$-central extension  $0 \to Z_{\Lie}(\Lieq) \to \Lieq \to \Lieq/Z_{\Lie}(\Lieq) \to 0$. \rdg

\begin{Rem}
If $\Lieq$ is a $\Lie$-unicentral Leibniz algebra, by Corollary \ref{charact} and the sequence $(\ref{six-term})$ associated to the $\Lie$-central extension  $0 \to Z_{\Lie}(\Lieq) \to \Lieq \to \Lieq/Z_{\Lie}(\Lieq) \to 0$, we have
\[
HL^{\Lie}_2(\Lieq) = {\rm Ker} \left( \theta^{\ast}(\Lieq) : HL^{\Lie}_2(\Lieq/Z_{\Lie}(\Lieq)) \to Z_{\Lie}(\Lieq) \right).
\]
\end{Rem}

\end{section}

%---------------------------------------------------------------------------------------

\begin{section}{$\Lie$-solvable and $\Lie$-nilpotent Leibniz algebras} \label{nil}

 In this section, by using the relative commutators, we introduce the notions of $\Lie$-solvability and $\Lie$-nilpotency of Leibniz algebras and investigate their properties.
\begin{De}
Let ${\Liem}$ be a two-sided ideal of a Leibniz algebra ${\Lieq}$. A series from ${\Liem}$ to ${\Lieq}$ is a finite sequence of two-sided ideals ${\Liem}_i$, $0 \leq i \leq k$, of ${\Lieq}$ such that
\[
{\Liem}= {\Liem}_0 \trianglelefteq {\Liem}_1 \trianglelefteq \dots \trianglelefteq {\Liem}_{k-1} \trianglelefteq
{\Liem}_k ={\Lieq} \; .
\]
$k$ is called the length of this series.
%The two-sided ideals ${\Liem}_i$ and the quotients
%${\Liem}_i/{\Liem}_{i-1}$ are called the terms and  the factors of the series, respectively.

   A series from $\Liem$ to $\Lieq$ of length $k$ is said to be $\Lie$-central (resp. $\Lie$-abelian) if $[{\Liem}_i, {\Lieq}]_{\Lie} \subseteq {\Liem}_{i-1}$,  or equivalently ${\Liem}_i/{\Liem}_{i-1} \subseteq Z_{\Lie}({\Lieq}/{\Liem}_{i-1})$ (resp. if $[\Liem_i,\Liem_i]_{\Lie} \subseteq \Liem_{i-1}$, or equivalently $[\Liem_i/\Liem_{i-1}, \Liem_i/\Liem_{i-1}]_{\Lie} = 0$) for $1\leq i\leq k$ .

A series from $0$ to ${\Lieq}$ is called a series of the Leibniz algebra ${\Lieq}$.
\end{De}

%\begin{Le} \label{ideal}
%Let $\Liem$ and $\Lien$ be  two-sided ideal of a Leibniz algebra \Lieq, then $[\Liem, \Lien]_{\Lie}$ is a two-sided ideal of $\Liem$, $\Lien$ and $\Lieq$.
%\end{Le}
%{\it Proof.} For any $m, m' \in \Liem, n \in \Lien$,
%$$[[m,n]+[n,m],m']=[m,[n,m']]+[[n,m'],m]+[[m,m'],n]+[n,[m,m']] \in [\Liem, \Lien]_{\Lie}$$
%and
%$$[m',[m,n]+[n,m]]=0.$$
%Other cases are completely analogous. \rdg

\begin{De}
A Leibniz algebra $\Lieq$ is said to be $\Lie$-solvable if it has a $\Lie$-abelian series.
If $k$ is the minimal length of such series, then $k$ is called the class of $\Lie$-solvability of ${\Lieq}$.
\end{De}

We show below that among all $\Lie$-abelian series of a $\Lie$-solvable Leibniz algebra there is one which descends most rapidly.

\begin{De}\label{central series}
The $\Lie$-derived series of a Leibniz algebra $\Lieq$  is the sequence
\[
\cdots \trianglelefteq {\Lieq}^{(i)} \trianglelefteq \cdots \trianglelefteq {\Lieq}^{(1)}  \trianglelefteq {\Lieq}^{(0)}
\]
of two-sided ideals of $\Lieq$ defined inductively by
 \[
 {\Lieq}^{(0)} = {\Lieq} \quad \text{and} \quad
 {\Lieq}^{(i)} =[{\Lieq}^{(i-1)},{\Lieq}^{(i-1)}]_{\Lie} ,  \   i \geq 1  .
 \]
\end{De}

\begin{Th} \
\begin{enumerate}
\item[(a)] Let $\Lieq$ be a Leibniz algebra and $\Liem = {\Liem}_0 \trianglelefteq {\Liem}_1 \trianglelefteq \dots \trianglelefteq {\Liem}_{j-1} \trianglelefteq {\Liem}_j ={\Lieq}$ be a $\Lie$-abelian series from $\Liem$ to ${\Lieq}$, then $ {\Lieq}^{(i)} \subseteq {\Liem}_{j-i}$, $0\leq i\leq j$.
\item[(b)] A Leibniz algebra ${\Lieq}$ is $\Lie$-solvable with class of $\Lie$-solvability $k$ if and only if ${\Lieq}^{(k)}=0$ and
${\Lieq}^{(k-1)}\ne 0$.
\end{enumerate}
\end{Th}
{\it Proof.} (a) This easily follows by induction on $i$.

(b)  If $\Lieq^{(k)} = 0$ and $\Lieq^{(k-1)} \neq 0$, then $0 = \Lieq^{(k)} \trianglelefteq \Lieq^{(k-1)} \trianglelefteq \dots \trianglelefteq \Lieq^{(1)} \trianglelefteq \Lieq^{(0)} ={\Lieq}$ is a $\Lie$-abelian series and by (a) its length is minimal. Therefore $\Lieq$ is a $\Lie$-solvable with class of solvability $k$. The converse statement directly follows from (a). \rdg

\begin{Ex}\label{ex solvable}\
\begin{enumerate}
\item[(a)] Abelian Leibniz algebras are  $\Lie$-solvable Leibniz algebras of class 1.
\item[(b)] Lie algebras are $\Lie$-solvable Leibniz algebras of class 1.
\item[(c)] The three-dimensional (non-Lie) Leibniz algebra with $\K$-linear basis $\{a_1,a_2,a_3 \}$ and Leibniz bracket given by $$[a_1,a_3]= [a_2,a_3]=a_2, [a_3,a_3]=a_1$$ and zero elsewhere (see \cite{CILL}), is a $\Lie$-solvable Leibniz algebra of class 2 and solvable of class 2.
\item[(d)] The five-dimensional perfect (non-Lie) Leibniz algebra with $\K$-linear basis $\{a_1,a_2,a_3, a_4, a_5 \}$ and Leibniz bracket given by
    \begin{center}
 $\begin{array}{lll}
 [a_2,a_1]=-a_3, & [a_1,a_2]= a_3,& [a_1,a_3]=-2a_1, \\

[a_3,a_1]= 2a_1, & [a_3,a_2]=-2a_2, & [a_2,a_3]=2a_2, \\

[a_5,a_1]=a_4, & [a_4,a_2]=a_5,& [a_4,a_3]=-a_4,\\

 & & [a_5,a_3]=a_5, \end{array}$
\end{center}

 and zero elsewhere (see \cite{Om}) is \Lie-solvable of class 2, but it is not a solvable Leibniz algebra.

\item[(e)] Subalgebras and images by homomorphisms of $\Lie$-solvable Leibniz algebras are $\Lie$-solvable as well.
\item[(f)] Intersection and sum of $\Lie$-solvable two-sided ideals of a Leibniz algebra are $\Lie$-solvable two-sided ideals as well.
\end{enumerate}
\end{Ex}

\begin{Pro}\
\begin{enumerate}
\item[(a)] Let $\Lien$ be a $\Lie$-solvable ideal of a Leibniz algebra $\Lieq$ such that $\Lieq / \Lien$ is $\Lie$-solvable, then $\Lieq$ itself is $\Lie$-solvable.

    \item[(b)] Let $\Liem$ and $\Lien$ be $\Lie$-solvable two-sided ideals of a Leibniz algebra $\Lieq$, then $[\Liem, \Lien]_{\Lie}$ is a $\Lie$-solvable two-sided ideal of $\Lieq$.
\end{enumerate}
\end{Pro}
{\it Proof.} (a) Since $\Lieq / \Lien$ is $\Lie$-solvable, then there exists $k \in \mathbb{N}$ such that $0 = \left( \Lieq / \Lien \right)^{(k)} = \Lieq^{(k)} / \Lien$, that is $ \Lieq^{(k)} =  \Lien$. On the other hand, since $\Lien$ is $\Lie$-solvable, then there exists $j \in \mathbb{N}$ such that $\Lien^{(j)}=0$. Hence $\Lieq^{(k+j)} = \left( \Lieq^{(k)} \right)^{(j)} = \Lien^{(j)}=0$.

(b) By Lemma \ref{Lemma3.3} $[\Liem, \Lien]_{\Lie}$ is a two-sided ideal of $\Lieq$. Moreover,  $[\Liem, \Lien]_{\Lie} \subseteq \Liem \cap \Lien$. Then the statements (e) and (f) in Example \ref{ex solvable} complete the proof. \rdg

\begin{De}
A Leibniz algebra $\Lieq$ is said to be $\Lie$-nilpotent if it has a $\Lie$-central series.
If $k$ is the minimal length of such series, then $k$ is
called the class of $\Lie$-nilpotency of ${\Lieq}$.
\end{De}

We show below that among all $\Lie$-central series of a $\Lie$-nilpotent Leibniz algebra there is one which descends most rapidly.

\begin{De}\label{lcs}
The lower $\Lie$-central series of a Leibniz algebra $\Lieq$  is the sequence
\[
\cdots \trianglelefteq{\Lieq}^{[i]} \trianglelefteq \cdots \trianglelefteq {\Lieq}^{[2]}  \trianglelefteq {\Lieq}^{[1]}
\]
of two-sided ideals of $\Lieq$ defined inductively by
 \[
 {\Lieq}^{[1]} = {\Lieq} \quad \text{and} \quad
{\Lieq}^{[i]} =[{\Lieq}^{[i-1]},{\Lieq}]_{\Lie} ,  \   i \geq 2  .
 \]
\end{De}

\begin{Th} \
\begin{enumerate}
\item[(a)] Let $\Lieq$ be a Leibniz algebra and $0={\Liem}_0 \trianglelefteq
{\Liem}_1 \trianglelefteq \dots \trianglelefteq{\Liem}_{j-1} \trianglelefteq {\Liem}_j ={\Lieq}$ be
a $\Lie$-central series of ${\Lieq}$, then $ {\Lieq}^{[i]} \subseteq
{\Liem}_{j-i+1}$, $1\leq i\leq j+1$.
\item[(b)] A Leibniz algebra ${\Lieq}$ is $\Lie$-nilpotent with class
of nilpotency $k$ if and only if ${\Lieq}^{[k+1]}=0$ and
${\Lieq}^{[k]}\ne 0$.
\end{enumerate}
\end{Th}
{\it Proof.} (a) This follows by induction on $i$.

(b) If ${\Lieq}^{[k+1]}=0$ and ${\Lieq}^{[k]}\ne 0$, then $0 ={\Lieq}^{[k+1]} \trianglelefteq {\Lieq}^{[k]} \trianglelefteq \dots \trianglelefteq {\Lieq}^{[2]} \trianglelefteq {\Lieq}^{[1]}= \Lieq$ is a $\Lie$-central series, which has minimal length by (a). Therefore $\Lieq$ is $\Lie$-nilpotent of class $k$. The inverse statement directly follows from (a). \rdg

\begin{De}
The upper $\Lie$-central series of a Leibniz algebra ${\Lieq}$ is the sequence of
 two-sided ideals
 \[
{\ze}_0^{\Lie}({\Lieq}) \trianglelefteq {\ze}_1^{\Lie}({\Lieq}) \trianglelefteq \cdots \trianglelefteq {\ze}_i^{\Lie}({\Lieq}) \trianglelefteq \cdots
\]
 defined inductively by
\[
{\ze}_0^{\Lie}({\Lieq}) = 0 \quad \text{and} \quad
 {\ze}_{i}^{\Lie}({\Lieq}) = C_{\Lieq}^{\Lie}({\Lieq},{\ze}_{i-1}^{\Lie}({\Lieq})) , \  i \geq 1 .
 \]
\end{De}

 Let us observe that ${\ze}_1^{\Lie}({\Lieq}) = Z_{\Lie}({\Lieq})$ and, by Lemma \ref{Lemma3.3}, ${\ze}_i^{\rm Lie}({\Lieq})$
  indeed is a two-sided ideal of ${\Lieq}$.

\begin{Le} \label{inclusion2}
Let ${\Liem}$ and ${\Lien}$ be two-sided ideals of a Leibniz algebra ${\Lieq}$.
 If $[{\Liem}, {\Lieq}]_{\Lie} \subseteq {\Lien}$, then ${\Liem}\subseteq C_{\Lieq}^{\Lie}({\Lieq},{\Lien})$.
\end{Le}
{\it Proof.} This is straightforward. \rdg
\bigskip

Now we show that among all $\Lie$-central series of a nilpotent Leibniz algebra there is one which ascends most rapidly.

\begin{Pro} \label{inclusion1}
Let $0 = {\Liem}_0 \trianglelefteq {\Liem}_1 \trianglelefteq \dots \trianglelefteq {\Liem}_{k-1} \trianglelefteq {\Liem}_k = {\Lieq}$ be a
$\Lie$-central series of a $\Lie$-nilpotent Leibniz algebra ${\Lieq}$, then
${\Liem}_i \subseteq {\ze}_i^{\Lie}({\Lieq})$, $0\leq i \leq k$.
\end{Pro}
{\it Proof.} Obviously the assertion holds for $i=0$. Proceeding by induction on $i$, we assume that the assertion is true for $i-1$, then using the $\Lie$-centrality  and Lemma \ref{inclusion2}, we have ${\Liem}_i \subseteq C_{\Lieq}^{\Lie}({\Lieq},{\Liem}_{i-1}) \subseteq C_{\Lieq}^{\Lie}({\Lieq},{\ze}_{i-1}^{\Lie}({\Lieq})) = {\ze}_{i}^{\Lie}({\Lieq}).$ \rdg

\begin{Th}
 A Leibniz algebra ${\Lieq}$ is $\Lie$-nilpotent with class of $\Lie$-nilpotency k if and only if
${\ze}_k^{\Lie}({\Lieq}) = {\Lieq}$ and ${\ze}_{k-1}^{\Lie}({\Lieq}) \ne {\Lieq}$.
\end{Th}
{\it Proof}. If ${\Lieq}$ is a $\Lie$-nilpotent Leibniz algebra with class of $\Lie$-nilpotency $k$, then  Proposition \ref{inclusion1} implies that  ${\Lieq} = {\Liem}_k \subseteq {\ze}_k^{\Lie}({\Lieq}) \subseteq {\Lieq}$. Moreover, in
this case $0 = {\ze}_0^{\Lie}({\Lieq}) \trianglelefteq {\ze}_1^{\Lie}({\Lieq}) \trianglelefteq \dots
\trianglelefteq {\ze}_{k-1}^{\Lie}({\Lieq}) \trianglelefteq {\ze}_k^{\Lie}({\Lieq}) = {\Lieq}$ is a $\Lie$-central series of length $k$ of ${\Lieq}$. Hence
${\ze}_{k-1}^{\Lie}({\Lieq}) \ne {\Lieq}$.

\noindent Conversely, if ${\ze}_k^{\Lie}({\Lieq}) = {\Lieq}$ and ${\ze}_{k-1}^{\Lie}({\Lieq}) \ne {\Lieq}$
then $0 = {\ze}_0^{\Lie}({\Lieq}) \trianglelefteq {\ze}_1^{\Lie}({\Lieq}) \trianglelefteq \dots \trianglelefteq {\ze}_{k-1}^{\Lie}({\Lieq}) \trianglelefteq {\ze}_k^{\Lie}({\Lieq}) = {\Lieq}$ is a $\Lie$-central series of ${\Lieq}$ and by Proposition \ref{inclusion1} its
length is minimal.  \rdg

\begin{Ex}\
\begin{enumerate}
\item[(a)] Abelian Leibniz algebras are  $\Lie$-nilpotent Leibniz algebras of class 1.
\item[(b)] Lie algebras are $\Lie$-nilpotent Leibniz algebras of class 1.
\item[(c)] The three-dimensional non-Lie Leibniz algebra  with $\K$-linear basis $\{a_1,a_2,a_3 \}$ and Leibniz bracket given by $$[a_3,a_3]=a_1$$ and zero elsewhere (see \cite{CILL}),  is a $\Lie$-nilpotent Leibniz algebra of class 2.
    \item[(d)] The non-Lie Leibniz algebra given in Example \ref{ex solvable} (c) is non $\Lie$-nilpotent.
 \item[(e)] The four-dimensional (non-Lie) Leibniz algebra with $\K$-linear basis $\{a_1,a_2,a_3, a_4 \}$ and Leibniz bracket given by $$[e_1,e_1]=e_3, [e_2,e_4]=e_2, [e_4,e_2]=-e_2$$  and zero elsewhere (see \cite{CK}) is \Lie-nilpotent of class 2, but it is not a nilpotent Leibniz algebra.
    \item[(f)]  Subalgebras and images by homomorphisms of $\Lie$-nilpotent Leibniz algebras are $\Lie$-nilpotent Leibniz algebras.
\item[(g)] Intersection and sum of $\Lie$-nilpotent two-sided ideals of a Leibniz algebra are $\Lie$-nilpotent two-sided ideals as well.
 \end{enumerate}
\end{Ex}

\begin{Pro}\
\begin{enumerate}
\item[(a)] If ${\Lieq}/Z_{\Lie}({\Lieq}$) is a $\Lie$-nilpotent Leibniz algebra, then ${\Lieq}$ is a $\Lie$-nilpotent Leibniz algebra.

\item[(b)] If ${\Lieq}$ is a $\Lie$-nilpotent and non trivial Leibniz algebra, then $Z_{\Lie}({\Lieq})\not = 0$.

\item[(c)] If $\Lieg\twoheadrightarrow\Lieq$ is a $\Lie$-central extension of a $\Lie$-nilpotent Leibniz algebra $\Lieq$, then $\Lieg$ is $\Lie$-nilpotent as well.

\item[(d)]  A $\Lie$-nilpotent Leibniz algebra is $\Lie$-solvable as well.
\end{enumerate}
\end{Pro}
{\it Proof.} (a) There exists $k \in \mathbb{N}$ such that $\left( \Lieq/Z_{\Lie}(\Lieq) \right)^{[k]}=0$, then $\Lieq^{[k]} \subseteq Z_{\Lie}(\Lieq)$, hence $\Lieq^{[k+1]} \subseteq [Z_{\Lie}(\Lieq), \Lieq]_{\Lie}=0$.

(b) Assume that ${\Lieq}$ has  $\Lie$-nilpotency class equal to $k$, that is $[{\Lieq}^{[k]},\Lieq]_{\Lie} ={\Lieq}^{[k+1]} = 0$, then Lemma \ref{inclusion2} implies that $0 \neq {\Lieq}^{[k]} \subseteq C_{\Lieq}^{\Lie}(\Lieq,0) = Z_{\Lie}(\Lieq)$.

(c) There exists $k \in \mathbb{N}$ such that $\Lieq^{[k]}=0$. Then $\left( \Lieg / \Lien \right)^{[k]} =\Lieg^{[k]}/ \Lien = 0$, where $\Lien = \Ker( \Lieg \twoheadrightarrow \Lieq)$. Hence $\Lieg^{[k]} \subseteq \Lien  \subseteq Z_{\Lie}(\Lieg)$ and $\Lieg^{[k+1]} = [\Lieg^{[k]},\Lieg]_{\Lie}=0$.

(d) By induction on $i$, it is easy to see that $\Lieq^{(i)} \subseteq \Lieq^{[i+1]}$, $i \geq 0$.
\rdg

%---------------------------------------------------------------------------------------

\section{Homological criterion for \Lie-nilpotency} \label{HC}

\begin{De}
Let ${\Lien}$ be a two-sided ideal of a Leibniz algebra $\Lieq$. The lower $\Lie$-central series of $\Lieq$ relative to ${\Lien}$ is the sequence
\[
\cdots \trianglelefteq {\Lien}^{[i]} \trianglelefteq \cdots \trianglelefteq {\Lien}^{[2]}  \trianglelefteq {\Lien}^{[1]}
\]
of two-sided ideals of $\Lieq$ defined inductively by
\[
{\Lien}^{[1]} = {\Lien} \quad \text{and} \quad {\Lien}^{[i]} =[{\Lien}^{[i-1]},{\Lieq}]_{\Lie}, \quad   i \geq 2.
\]

\end{De}

Note that $[{\Lien}^{[i]}/{\Lien}^{[i+1]},{\Lien}^{[i]}/{\Lien}^{[i+1]}]_{\Lie}={0}$. When ${\Lien} = {\Lieq}$ we obtain  Definition \ref{lcs}.
If $\varphi : {\Lieg} \to {\Lieq}$ is a homomorphism of Leibniz algebras such that $\varphi({\Liem}) \subseteq {\Lien}$, where
${\Liem}$ is a two-sided ideal of ${\Lieg}$ and ${\Lien}$ is a two-sided ideal of ${\Lieq}$, then $\varphi({\Liem}^{[i]}) \subseteq
{\Lien}^{[i]}$, $i \geq 1$.

\begin{Th} \label{teorema five-term}
Let $\varphi : {\Lieg} \to {\Lieq}$ be a homomorphism of Leibniz algebras such that $\varphi({\Liem}) \subseteq {\Lien}$, where ${\Liem}$ is a two-sided ideal of ${\Lieg}$ and ${\Lien}$ is a two-sided ideal of ${\Lieq}$, and
 the following properties hold:
\begin{enumerate}
\item[(a)] the induced homomorphism $HL^{\Lie}_1({\Lieg}) \to HL^{\Lie}_1({\Lieq})$ is an isomorphism;
\item[(b)]  the induced homomorphism $HL^{\Lie}_2({\Lieg}) \to HL^{\Lie}_2({\Lieq})$ is an epimorphism;
\item[(c)] the induced homomorphism ${\varphi}_1 : {\Lieg}/{\Liem} \to {\Lieq}/{\Lien}$ is an isomorphism.
\end{enumerate}
Then $\varphi$ induces a natural isomorphism $\varphi_k : {\Lieg}/{\Liem}^{[k]} \to {\Lieq}/{\Lien}^{[k]}$, $k \geq 1$.
\end{Th}
{\it Proof.} We prove by induction on $k$. For $k=1$ we have the statement (c). Suppose  the theorem is true for $k-1$. Applying
sequence (\ref{five-term trivial}), which is natural \cite{EVDL},  to the following commutative diagram
\[ \xymatrix{
0  \ar[r] & {\Liem}^{[k-1]} \ar[r] \ar[d] & {\Lieg} \ar[r] \ar[d]^{\varphi} & {\Lieg}/{\Liem}^{[k-1]} \ar[r] \ar[d]_{\wr}^{\varphi_{k-1}} & 0 \\
0  \ar[r] & {\Lien}^{[k-1]} \ar[r]  & {\Lieq} \ar[r]  & {\Lieq}/{\Lien}^{[k-1]} \ar[r]  & 0
} \]
 we get the following commutative diagram
\[ \xymatrix{
HL^{\Lie}_2(\Lieg) \ar[r] \ar[d]& HL^{\Lie}_2(\frac{\Lieg}{\Liem^{[k-1]}}) \ar[r] \ar[d] & \frac{\Liem^{[k-1]}}{[\Liem^{[k-1]},\Lieg]_{\Lie}} \ar[r] \ar[d] & HL^{\Lie}_1(\Lieg) \ar[r] \ar[d]& HL^{\Lie}_1(\frac{\Lieg}{\Liem^{[k-1]}}) \ar[r] \ar[d] & 0\\
HL^{\Lie}_2(\Lieq) \ar[r] & HL^{\Lie}_2(\frac{\Lieq}{\Lien^{[k-1]}}) \ar[r]  & \frac{\Lien^{[k-1]}}{[\Lien^{[k-1]},\Lieq]_{\Lie}} \ar[r]  & HL^{\Lie}_1(\Lieq) \ar[r] & HL^{\Lie}_1(\frac{\Lieq}{\Lien^{[k-1]}}) \ar[r]  & 0
} \]
By the Five Lemma, which holds in a semi-abelian category  \cite{Bor, Mi}, we get
 $\frac{\Liem^{[k-1]}}{[\Liem^{[k-1]},\Lieg]_{\Lie}} \cong \frac{\Lien^{[k-1]}}{[\Lien^{[k-1]},\Lieq]_{\Lie}}$, i. e.
 $\frac{\Liem^{[k-1]}}{\Liem^{[k]}} \cong \frac{\Lien^{[k-1]}}{\Lien^{[k]}}$.
Then the short Five Lemma applied to the  following commutative diagram
\[ \xymatrix{
0 \ar[r] & \frac{\Liem^{[k-1]}}{\Liem^{[k]}} \ar[r] \ar[d] & \frac{\Lieg}{\Liem^{[k]}} \ar[r] \ar[d] & \frac{\Lieg}{\Liem^{[k-1]}} \ar[r] \ar[d] & 0 \\
0 \ar[r] & \frac{\Lien^{[k-1]}}{\Lien^{[k]}} \ar[r] & \frac{\Lieq}{\Lien^{[k]}} \ar[r]
 & \frac{\Lieq}{\Lien^{[k-1]}} \ar[r]  & 0
} \]
and the induction complete the proof. \rdg

\begin{Co} \label{isom}
Let $\varphi : {\Lieg} \to {\Lieq}$ be a homomorphism of Leibniz algebras such that $\varphi_{\rm Lie} : {\Lieg}_{\rm Lie} \to {\Lieq}_{\rm Lie}$ is an isomorphism and the induced homomorphism $HL^{\Lie}_2({\Lieg}) \to HL^{\Lie}_2({\Lieq})$ is an epimorphism.
If ${\Lieg}$ and ${\Lieq}$ are $\Lie$-nilpotent Leibniz algebras, then $\varphi$ is an isomorphism.
\end{Co}
{\it Proof.} Take ${\Liem}={\Lieg}$ and ${\Lien} = {\Lieq}$ in Theorem \ref{teorema five-term}. Then the assertion follows by keeping in mind that $HL^{\Lie}_1({\Lieg}) \cong {\Lieg}_{\rm Lie}$, $HL^{\Lie}_1({\Lieq}) \cong {\Lieq}_{\rm Lie}$  and there exists $k\geq 1$ such that ${\Lieg}^{[k]} = {\Lieq}^{[k]} = 0$. \rdg

\begin{Co}
Suppose ${\Lieq}$ and ${\Lieq}_{\rm Lie}$ both are $\Lie$-nilpotent Leibniz algebras. Then necessarily  ${\Lieq}
\cong {\Lieq}_{\rm Lie}$, that is,  $\Lieq$ is a Lie algebra.
\end{Co}
{\it Proof.} This follows by applying Corollary \ref{isom} to the
canonical epimorphism ${\Lieq} \twoheadrightarrow {\Lieq}_{\rm Lie}$. \rdg

\begin{Le} \label{inclusion}
Let ${\Lien}$ be a two-sided ideal of a Leibniz algebra ${\Lieq}$. The lower $\Lie$-central series determined by ${\Lien}$ vanishes
if and only if there exists $i \geq 0$ such that ${\Lien} \subseteq {\ze}_{i}^{\Lie}({\Lieq})$.
\end{Le}
{\it Proof.} It is enough to use the following  obvious
equivalence:
$${\Lien}^{[i]} \subseteq {\ze}_k^{\Lie}({\Lieq})  \Leftrightarrow
{\Lien}^{[i+1]} \subseteq {\ze}_{k-1}^{\Lie}({\Lieq}).$$ \rdg

\begin{Th} \label{characterization}
Let $\Liem$ be a two-sided ideal of $\Lieg$ and $\Lien$ be a two-sided ideal of $\Lieq$ such that  ${\Liem} \subseteq {\ze}_i^{\Lie}({\Lieg})$ and ${\Lien} \subseteq {\ze}_j^{\Lie}({\Lieq})$, $i,j \geq 0$. Let $\varphi : {\Lieg} \to {\Lieq}$ be a homomorphism of Leibniz algebras such that $\varphi({\Liem}) \subseteq {\Lien}$ and  satisfying the following conditions:
\begin{enumerate}
\item[(a)] the induced homomorphism $HL^{\Lie}_1({\Lieg}) \to HL^{\Lie}_1({\Lieq})$ is an isomorphism;

\item[(b)] the induced homomorphism $HL^{\Lie}_2({\Lieg}) \to HL^{\Lie}_2({\Lieq})$ is an epimorphism;

\item[(c)] the induced homomorphism ${\Lieg}/{\Liem} \to {\Lieq}/{\Lien}$ is an isomorphism.
\end{enumerate}
Then $\varphi : {\Lieg} \to {\Lieq}$ is an isomorphism.
\end{Th}
{\it Proof.} This is a consequence of Theorem \ref{teorema five-term} and Lemma \ref{inclusion}. \rdg

\begin{Co}
 Let $\varphi : {\Lieg} \to {\Lieq}$ be a homomorphism of Leibniz algebras such that $\varphi({\ze}_i^{\Lie}({\Lieg})) \subseteq
{\ze}_i^{\Lie}({\Lieq})$ for any $i \geq 0$ and satisfying the following conditions:
\begin{enumerate}
\item[(a)] the induced homomorphism $HL^{\Lie}_1({\Lieg}) \to HL^{\Lie}_1({\Lieq})$ is an isomorphism;

\item[(b)] the induced homomorphism $HL^{\Lie}_2({\Lieg}) \to HL^{\Lie}_2({\Lieq})$ is an epimorphism;

\item[(c)] the induced homomorphism ${\Lieg}/{\ze}_i^{\Lie}({\Lieg}) \to {\Lieq}/{\ze}_i^{\Lie}({\Lieq})$ is an
isomorphism.
\end{enumerate}
Then $\varphi : {\Lieg} \to {\Lieq}$ is an isomorphism.
\end{Co}
{\it Proof.} This follows by applying Theorem \ref{characterization} to the case ${\Liem}={\ze}_i^{\Lie}({\Lieg})$
 and ${\Lien} = {\ze}_i^{\Lie}({\Lieq})$. \rdg

\end{section}
\bigskip
%---------------------------------------------------------------------------------------

\centerline{\bf Acknowledgements}
The authors wish to thank Prof. Teiumuraz Pirashvili for his helpful suggestions which have improved the content of the paper.
The authors were supported by Ministerio de Economía y Competitividad (Spain) (European FEDER support included), grant MTM2013-43687-P.  The second author was supported by  Shota Rustaveli National Science Foundation, grant FR/189/5-113/14.

%---------------------------------------------------------------------------------------

\begin{center}

\end{center}

\end{document}